\documentclass[12pt,a4paper]{amsart}
\usepackage{amssymb}
\usepackage{mathrsfs}
\headsep=1cm \numberwithin{equation}{section}
\newtheorem{thm}{Theorem}[section]
\newtheorem{cor}[thm]{Corollary}
\newtheorem{lem}[thm]{Lemma}

\newtheorem{prop}[thm]{Proposition}
\theoremstyle{definition}
\newtheorem{defn}[thm]{Definition}

\newtheorem{rem}[thm]{Remark}
\numberwithin{equation}{section}

\newcommand\Supp{\operatorname{Supp}}
\newcommand\mAtt{\operatorname{mAtt}}
\newcommand\cd{\operatorname{cd}}

\newcommand\ara{\operatorname{ara}}

\newcommand\Ass{\operatorname{Ass}}
\newcommand\mAss{\operatorname{mAss}}

\newcommand\Ann{\operatorname{Ann}}
\newcommand\Spec{\operatorname{Spec}}
\newcommand\Rad{\operatorname{Rad}}

\newcommand\Ext{\operatorname{Ext}}
\newcommand\grade{\operatorname{grade}}

\newcommand\Att{\operatorname{Att}}

\newcommand\V{\operatorname{V}}

\begin{document}\title[Annihilators and attached primes of local cohomology]
{Some results on the annihilators and attached primes of local cohomology modules}
\author{Ali Atazadeh, Monireh Sedghi$^*$ and Reza Naghipour
\vspace*{0.2cm}}
\address{Department of Mathematics, Kaleybar Branch, Islamic Azad University, Kaleybar, Iran; and
Department of Mathematics, Azarbaijan Shahid Madani University, Tabriz, Iran.} 
\email{aalzp2002@yahoo.com}

\address{Department of Mathematics, Azarbaijan Shahid Madani University, Tabriz, Iran.}
\email {m\_sedghi@tabrizu.ac.ir}
\email {sedghi@azaruniv.ac.ir}

\address{Department of Mathematics, University of Tabriz, Tabriz, Iran.}
\email{naghipour@ipm.ir} \email {naghipour@tabrizu.ac.ir}
\thanks{ 
This research  was in part supported by the Azarbaijan Shahid Madani University.\\
$^*$Corresponding author: e-mail: {\it m\_sedghi@tabrizu.ac.ir and sedghi@azaruniv.ac.ir} (Monireh Sedghi)}%

\begin{abstract}
Let $(R, \frak m)$ be a local  ring  and $M$ a finitely generated $R$-module.  It is shown that if  $M$ is  relative Cohen-Macaulay  with respect to an ideal $\frak a$ of $R$,
then $\Ann_R(H_{\frak a}^{\cd(\frak a, M)}(M))=\Ann_RM/L=\Ann_RM$ and $\Ass_R (R/\Ann_RM)\subseteq \{\frak p\in \Ass_R M|\,{\rm cd}(\frak a, R/\frak p)=\cd(\frak a, M)\},$
where  $L$ is the largest submodule of $M$ such that ${\rm cd}(\frak a, L)< {\rm cd}(\frak a, M)$. We also show that if  $H^{\dim M}_{\frak a}(M)=0$, then
$\Att_R(H^{\dim M-1}_{\frak a}(M))= \{\frak p\in \Supp (M)|\,{\rm cd}(\frak a, R/\frak p)=\dim M-1\},$ and so the attached primes of  $H^{\dim M-1}_{\frak a}(M)$
 depends only on $\Supp (M)$.  Finally, we prove that if  $M$  is an arbitrary module (not necessarily finitely generated) over a Noetherian ring $R$ with
  ${\rm cd}(\frak a, M)={\rm cd}(\frak a, R/\Ann_RM)$, then
 $\Att_R(H^{{\rm cd}(\frak a, M)}_{\frak a}(M))\subseteq\{\frak p\in \V(\Ann_RM)|\,{\rm cd}(\frak a, R/\frak p)={\rm cd}(\frak a, M)\}.$
  As a consequence of this it is shown that if $\dim M=\dim R$, then
 $\Att_R(H^{\dim M}_{\frak a}(M))\subseteq\{\frak p\in \Ass_R M|\,{\rm cd}(\frak a, R/\frak p)=\dim M\}.$

\end{abstract}
\maketitle
\section{Introduction}
Let $R$ be an arbitrary commutative Noetherian ring (with identity), $\frak a$ an ideal of $R$ and let $M$ be a finitely generated $R$-module.
Recall that  the $i^{\rm th}$ local cohomology module of $M$ with support in $V(\frak a)$ is defined by
$H^i_{\frak a}(M) := \underset{n\geq1} {\varinjlim}\,\, \Ext^i_R(R/\frak a^n, M).$ As the structure of local cohomology modules in general seems to be quite
mysterious, one tries to establish some properties providing a better understanding of these modules. Among these properties, an interesting question  is
 determining the annihilators of  local cohomology modules. This problem has been studied by several authors; see for example
 \cite{ASN1},  \cite{ASN2}, \cite{HK}, \cite{Lyn}, \cite{Ly1}, \cite{NA},  \cite{NC},  \cite{Sc1},   and  has led to some interesting results.
   A very interesting result shows that if $R$ is regular local ring containing a field, then $H^i_{\frak a}(R)\neq 0$, if and only if
    $\Ann_R(H_{\frak a}^i(R))=0$, for all $i\geq 0$,  cf. \cite{HK} (in positive characteristic) and \cite{Ly1} (in characteristic zero).
    One  purpose of the present paper  is to establish some new results concerning of the annihilators of local cohomology modules
    $H_{\frak a}^i(M)$ $(i\in \mathbb{N}_0)$. As a main result in the second section, we determine the  annihilators of the local cohomology module
     $H_{\frak a}^i(M)$ in several cases.  More precisely, we shall prove the following theorem:
\begin{thm}
Let $R$ be a Noetherian ring, $\frak a$ an ideal of $R$ and $M$  a non-zero finitely generated $R$-module  such that ${\rm cd}(\frak a, M)=c$.

$\rm(a)$ If $\dim M=c$, then
\begin{itemize}
  \item[(i)]
   $\Ann_R(H_{\frak a}^{c}(M))=\cap_{{\rm cd}(\frak a, R/{\frak p_i})=c}\,\frak q_i$, where $\Ann_RM=\cap_{i=1}^m\frak q_i$ denotes a
   reduced primary decomposition of the ideal $\Ann_RM$ and $\frak q_i$ is a $\frak p_i$-primary ideal of $R$, for all  $i=1,\dots ,m$.

\item[(ii)] $\Ann_R(H_{\frak a}^c(M))=\Ann_RM$ if and only if $\Ass_R (R/\Ann_RM) =\Att_RH^c_{\frak a}(M)$.
\end{itemize}

$\rm(b)$  If $(R, \frak m)$  is local and  $M$ is  relative Cohen-Macaulay with respect to  $\frak a$, then
\begin{itemize}
  \item[(i)]
$\Ann_R(H_{\frak a}^c(M))=\Ann_RM/T_R(\frak a, M)=\Ann_RM,$
where $T_R(\frak a, M)$ denotes   the largest submodule of $M$ such that ${\rm cd}(\frak a, T_R(\frak a, M))<{\rm cd}(\frak a, M)$.

\item[(ii)]
$\Ass_R (R/\Ann_RM)\subseteq \{\frak p\in \Ass_RM|\,{\rm cd}(\frak a, R/\frak p)=c\}$.
\end{itemize}
\end{thm}

The result  in Theorem 1.1 is proved in Theorems 2.2,  2.11 and Corollaries 2.3, 2.17.  One of our tools for proving Theorem 1.1 is the following, which plays a key role in this section.
\begin{prop}
Let $R$ be a local (Noetherian) ring and let ${\bf x}=x_1,\dots,x_n$ be elements of $R$.  Then for any  finitely generated
$R$-module $M$, $$(\Ann_R M)^{n+1}\subseteq (\Ann_R(H_{{\bf x}R}^0(M)))\cdots (\Ann_R(H_{{\bf x}R}^n(M)))\subseteq \Ann_RM.$$
\end{prop}

Another basic problems concerning local cohomology is to find the set of attached primes of  $H_{\frak a}^{i}(M)$.
In  Section 3, we obtain some results about the attached primes of  local cohomology modules. In this section among other things, we derive the following result,
  which is a generalization of the main results of \cite{DY} and \cite{Ru} in the case that $\dim M=\dim R$, to arbitrary modules $M$ (which may not be finitely generated)
  over an arbitrary Noetherian (not necessarily local)  ring $R$.
Recall that the dimension $\dim M$ of $M$ is the supremum of lengths of chains of prime ideals in $\Supp(M)$, if  this supremum exists, and $\infty$ otherwise.
\begin{thm}
Let $R$ be an arbitrary Noetherian ring and $\frak a$ an ideal of $R$. Let $M$ be an $R$-module (not necessarily finitely generated) with $\dim M=\dim R$. Then
 $$\Att_R(H^{\dim M}_{\frak a}(M))\subseteq\{\frak p\in \Ass_R M|\,{\rm cd}(\frak a, R/\frak p)=\dim M\}.$$
\end{thm}
The result  in Theorem 1.3 is proved in  Corollary 3.14.  One of our tools for proving Theorem 1.3 is the following.
\begin{prop}
Let $R$ be a Noetherian ring and $\frak a$ an ideal of $R$. Let $M$ be an $R$-module (not necessarily finitely generated)  such that ${\rm cd}(\frak a, M)={\rm cd}(\frak a, R/\Ann_RM)$.
Then
 $$\Att_R(H^{{\rm cd}(\frak a, M)}_{\frak a}(M))\subseteq\{\frak p\in \V(\Ann_RM)|\,{\rm cd}(\frak a, R/\frak p)={\rm cd}(\frak a, M)\}.$$
\end{prop}

For  an  $R$-module $A$, a prime ideal $\frak p$ of $R$  is said to be {\it attached prime to} $A$ if $\frak p=\Ann_R(A/B)$ for some submodule $B$
of $A$. We denote the set of attached primes of $A$ by $\Att_RA$.  This definition agrees with the usual definition of  attached  prime if $A$ has a
secondary representation (cf. \cite[Theorem 2.5]{Ma}).

Another main result in Section 3 is to give a complete characterization of  the attached primes of the local cohomology module $H_{\frak a}^{\dim M-1}(M)$.
More precisely, we shall show the following result,  which is an extension of the main theorems of \cite{ASN2} and \cite{He}.

\begin{thm}
 Let $R$ be a Noetherian ring  and $\frak a$ an ideal of $R$. Let  $M$ be a finitely generated $R$-module such that $H^{\dim M}_{\frak a}(M)=0$. Then

$$\Att_R(H^{\dim M-1}_{\frak a}(M))= \{\frak p\in \Supp (M)|\,{\rm cd}(\frak a, R/\frak p)=\dim M-1\}.$$
\end{thm}
As a consequence of Theorem 1.5 we show that the set $\Att_R(H_{\frak a}^{{\dim M-1}}(M))$ depends on $\Supp(M)$ only, whenever $H_{\frak a}^{\dim M}(M)=0.$
 More precisely, we shall show that:

\begin{cor}
Let $R$ be a Noetherian ring and $\frak a$ an ideal of $R$. Let $M$ and $N$ be two non-zero finitely generated $R$-modules with
$\dim M=d$ and $H_{\frak a}^d(M)=0.$ If $\,\Supp (M)=\Supp (N)$, then $\Att_R(H_{\frak a}^{d-1}(N))=\Att_R(H_{\frak a}^{d-1}(M)).$
\end{cor}

Throughout this paper, $R$ will always be a commutative Noetherian
ring with non-zero identity and $\frak a$ will be an ideal of $R$.
For each $R$-module $L$, we denote by
 ${\rm Assh}_RL$ (resp. $\mAss_RL$) the set $\{\frak p\in \Ass
_RL:\, \dim R/\frak p= \dim L\}$ (resp. the set of minimal primes of
$\Ass_RL$).  Also,  we shall use $\Att_R L$ to denote the set of attached prime
ideals of $L$.
For any ideal $\frak a$ of $R$, we denote $\{\frak p \in {\rm Spec}\,R:\, \frak p\supseteq \frak a \}$ by
$V(\frak a)$. Finally, for any ideal $\frak{b}$ of $R$, {\it the radical} of $\frak{b}$, denoted by $\Rad(\frak{b})$, is defined to
be the set $\{x\in R \,: \, x^n \in \frak{b}$ for some $n \in \mathbb{N}\}$.  For any unexplained notation and terminology we refer
the reader to \cite{BS} and \cite{Mat}.\\

\section{Annihilators of top local cohomology modules}

Let us, firstly, recall the important concept of the cohomological dimension of an $R$-module $L$ with respect to an ideal $\frak a$ of
a commutative Noetherian ring $R$, denoted by ${\rm cd}(\frak a,L)$, is the largest integer $i$ such that $H^i_{\frak a}(L)\neq 0$; i.e.,
 ${\rm cd}(\frak a, L):= \sup\{i\in \mathbb{Z}|\, H^i_{\frak a}(L)\neq 0\}.$  The first main observation of this section is Theorem {\rm2.2}.
  The following lemma plays a key role in the proof of that theorem.

\begin{lem}
Let $R$ be a Noetherian ring,  $\frak a$  an ideal of $R$ and let $M$ be a finitely generated $R$-module  with finite dimension $c$ such that ${\rm cd}(\frak a, M)=c$. Then
 $$\Ann_R(H_{\frak a}^c(M))=\Ann_R(M/H_{\frak b}^{0}(M))=\Ann_R(M/\cap_{{\rm cd}(\frak a, R/{\frak p_j})=c}N_j).$$
Here $0=\cap_{j=1}^nN_j$ denotes a reduced primary decomposition of zero submodule $0$ in $M$ and $N_j$ is a $\frak p_j$-primary submodule of $M$,
for all $j=1,\dots ,n$ and $\frak b=\Pi_{{\rm cd}(\frak a, R/{\frak p_j})\neq c}\,\frak p_j.$
\end{lem}
\proof   See  \cite[Corollary 2.7]{ASN2}. \qed\\

We are now ready to state and prove the first main result of this section.
\begin{thm}
Let $R$ be a Noetherian ring and $\frak a$ an ideal of $R$. Let $M$ be a non-zero finitely generated $R$-module with finite dimension $c$
such that ${\rm cd}(\frak a, M)=c$. Then $$\Ann_R(H_{\frak a}^{c}(M))=\cap_{{\rm cd}(\frak a, R/{\frak p_i})=c}\,\frak q_i.$$ Here $\Ann_RM=\cap_{i=1}^m\frak q_i$
 denotes a reduced primary decomposition of the ideal $\Ann_RM$ and $\frak q_i$ is a $\frak p_i$-primary ideal of $R$, for all  $i=1,\dots ,m.$
\end{thm}
\proof Let $0=\cap_{i=1}^nN_i$ denote a reduced primary decomposition of zero submodule $0$ in $M$ and $N_i$ is a $\frak p'_i$-primary submodule of $M$, for all $i=1,\dots ,n$.
Since $\mAss_RM=\mAss_R(R/\Ann_RM)$, without loss of generality we may (and do) assume that, there is some $t\leq \min \{m,\,n\}$,
such that $\mAss_RM=\{\frak p _1',\dots, \frak p' _t\}$, $\mAss_R(R/\Ann_RM)=\{\frak p _1,\dots, \frak p _t\}$ and $\frak p' _i=\frak p_i$, for each $i\leq t$.
On the other hand, it follows from $0=\cap_{i=1}^nN_i$ that  $\Ann_RM=\bigcap_{i=1}^n(N_i:_RM)$ is a primary decomposition (not necessarily reduced) of $\Ann_RM$,
 and so  $\frak q _i=(N_i:_RM)$, for each $i\leq t$. Now, in view of Lemma {\rm2.1} we have
\begin{eqnarray*}
\Ann_R(H_{\frak a}^{c}(M)) &=&\Ann_R(M/\cap_{\frak p'_i\in \Ass_R M,\,{\rm cd}(\frak a, R/{\frak p'_i})=c}N_i)
\\&=& \cap_{\frak p'_i\in \Ass_R M,\,{\rm cd}(\frak a, R/{\frak p'_i})=c}(N_i:_RM) \\&=& \cap_{ \frak p_i\in \Ass_R (R/\Ann_RM),\,{\rm cd}(\frak a, R/{\frak p_i})=c}\,\frak q_i,
\end{eqnarray*}
 as required.\qed\\

The first application of Theorem 2.2 improves  the main results of  \cite[Corollary 2.9]{BAG} and \cite[Corollary 2.5]{Lyn}.
 \begin{cor}
Let $R$ be a Noetherian ring and $\frak a$ an ideal of $R$. Let $M$ be a non-zero finitely generated $R$-module with finite dimension $c$ such that ${\rm cd}(\frak a, M)=c$.
Then the following conditions are equivalent:

$\rm(i)$ $\Ann_RH_{\frak a}^c(M)=\Ann_RM$.

$\rm(ii)$ $\Ass_R (R/\Ann_RM) =\Att_RH^c_{\frak a}(M)$.\\
In particular, when $M$ is a faithful $R$-module we have $\Ann_RH_{\frak a}^c(M)=0$ if and only if $\Ass_RR=\Att_RH^c_{\frak a}(M).$
\end{cor}
\proof The assertion follows from Theorem {\rm2.2} and \cite[Corollary 3.4]{ASN2}.  \qed\\

The following lemma will be used in  Proposition 2.5 and Theorem 2.6.
\begin{lem}
Let $R$ be a Noetherian ring and $\frak a$ an ideal of $R$. Let $M$ and $N$ be two non-zero finitely generated
$R$-modules such that $\Supp (N)\subseteq \Supp (M)$. Then  ${\rm cd}(\frak a, N)\leq {\rm cd}(\frak a, M).$
\end{lem}
\proof See \cite[Theorem 2.2]{DNT}. \qed\\

\begin{prop}
Let $R$ be a Noetherian ring and let  $\frak a, \frak b, \frak c$ be ideals of $R$. Let $M$ be a finitely generated
$R$-module such that $H^c_{\frak a}(M/\frak a M)\cong H^c_{\frak a}(M/\frak b M)\cong H^c_{\frak a}(M),$ where  $c:=\cd(\frak a, M)$ is finite.
 Then $$H^c_{\frak a}(M)\cong H^c_{\frak a}(M/(\frak a M+\frak b M)).$$
\end{prop}
\proof Since $$H^c_{\frak a}(M/\frak a M)\cong H^c_{\frak a}(M/\frak b M)\cong H^c_{\frak a}(M),$$  it yields that  $(\frak a+\frak b)H^c_{\frak a}(M)=0$.
Now as, $\Supp (M)=\Supp (R/\Ann_RM)$, it follows from Lemma {\rm2.4} that $\cd(\frak a, R/\Ann_RM)=c.$ Now it follows from \cite[Exercise 24]{Hu} and
 Independence theorem (cf.  \cite[Theorem 4.2.1]{BS}) that
\begin{eqnarray*}
H^c_{\frak a}(M/(\frak b M+\frak c M)) &\cong&H^c_{\frak a(R/\Ann_RM)}(R/\Ann_RM\otimes_{R/\Ann_RM}M/(\frak b M+\frak c M))
\\&\cong& H^c_{\frak a(R/\Ann_RM)}(R/\Ann_RM)\otimes_{R/\Ann_RM}M/(\frak b M+\frak c M)
\\&\cong& H^c_{\frak a(R/\Ann_RM)}(R/\Ann_RM)\otimes_{R/\Ann_RM}M\otimes_RR/(\frak b+\frak c)
\\&\cong&H^c_{\frak a(R/\Ann_RM)}(R/\Ann_RM\otimes_{R/\Ann_RM}M)\otimes_RR/(\frak b+\frak c)
\\&\cong&H^c_{\frak a}(M)\otimes_R(R/(\frak b+\frak c))\cong H^c_{\frak a}(M)/(\frak b+\frak c)H^c_{\frak a}(M)=H^c_{\frak a}(M),
\end{eqnarray*}
as required.\qed\\

\begin{thm}
Let $R$ be a Noetherian ring and $M$   a finitely generated
$R$-module. Let $\frak a$ be an ideal of $R$ such that $c:=\cd(\frak a, M)$ is finite. Then the set
$$\Sigma:=\{\frak c|  \frak c \text{ is an ideal of R and}\,\,  H^c_{\frak a}(M)\cong H^c_{\frak a}(M/\frak c M)\}$$
has a largest  element with respect to inclusion,  $\frak b$ say; and  $\frak b=\Ann_R(H^c_{\frak a}(M))$.
\end{thm}
\proof Since  $R$ is a Noetherian ring, the set $\Sigma$ has a maximal member, $\frak b$ say. (Note that  $\Ann_RM\in \Sigma$, and  so $\Sigma$  is not empty.)
 Since by Proposition {\rm2.5} the sum of any two members of $\Sigma$ is again in $\Sigma$, it follows that $\frak b$ contains every member of $\Sigma$,
 and so is the largest  element of $\Sigma$. Now, we show that $\frak b=\Ann_R(H^c_{\frak a}(M))$. To this end, as
$ H^c_{\frak a}(M)\cong H^c_{\frak a}(M/\frak b M)$, it follows that $\frak b\subseteq\Ann_R(H^c_{\frak a}(M))$. To establish the reverse inclusion, let
 $x\in \Ann_R(H_{\frak a}^c(M))$, and we show that $x\in \frak b$.  Because of $\frak b$ is the largest element of $\Sigma$,
 it is enough for us to show that $H^c_{\frak a}(M)\cong H^c_{\frak a}(M/xM)$.  As $\Supp (M)=\Supp (R/\Ann_RM)$,
  it follows from Lemma {\rm2.4} that $\cd(\frak a, R/\Ann_RM)=c$, and so it follows from \cite[Exercise 24]{Hu} and Independence theorem (see \cite[Theorem 4.2.1]{BS}) that
\begin{eqnarray*}
H^c_{\frak a}(M/xM) &\cong&H^c_{\frak a(R/\Ann_RM)}(R/\Ann_RM\otimes_{R/\Ann_RM}M/xM)
\\&\cong& H^c_{\frak a(R/\Ann_RM)}(R/\Ann_RM)\otimes_{R/\Ann_RM}M/xM \\&\cong& H^c_{\frak a(R/\Ann_RM)}(R/\Ann_RM)\otimes_{R/\Ann_RM}M\otimes_RR/xR
\\&\cong&H^c_{\frak a(R/\Ann_RM)}(R/\Ann_RM\otimes_{R/\Ann_RM}M)\otimes_RR/xR\\&\cong&H^c_{\frak a}(M)\otimes_R R/xR\cong H^c_{\frak a}(M)/xH^c_{\frak a}(M)=H^c_{\frak a}(M),
\end{eqnarray*}
as required.\qed\\

The following result  follows by the similar argument as in the proof of  \cite[Lemma 2.4.4]{Sc2}, but we give a direct proof for the convenience of the reader.

\begin{prop}
Let $R$ be a local (Noetherian) ring and let ${\bf x}=x_1,\dots,x_n$ be elements of $R$. Let $\frak a$ be an ideal of $R$ such that $\Rad(\frak a)=\Rad({\bf x}R)$.
 Then for any  finitely generated
$R$-module $M$, $$(\Ann_R M)^{n+1}\subseteq (\Ann_R(H_{\frak a}^0(M)))\cdots (\Ann_R(H_{\frak a}^n(M)))\subseteq \Ann_RM.$$
\end{prop}
\proof Without loss of generality we may assume that $\frak a={\bf
x}R$. The first containment is trivial, note that $\Ann_RM\subseteq
\Ann_R(H_{\frak a}^i(M))$, for all $i$. For the second one, let
$K^{\centerdot}({\bf x}^t, M)$ denote the Koszul complex of $M$ with
respect to ${\bf x}^t=x_1^t,\dots,x_n^t$, for each positive integer
$t$. Because of the cohomology modules $H^i({\bf x}^t, M)\, ,i\in
\mathbb{Z}$, of the Koszul complex $K^{\centerdot}({\bf x}^t, M)$
are annihilated by ${\bf x}^tR$, it follows that the support of
$H^i({\bf x}^t, M)\, ,i\in \mathbb{Z}$, is contained in $V(\frak
a)$. Now, in view of \cite[Corollary 1]{Sc1},
   $$(\Ann_R(H_{\frak a}^0(M)))\cdots (\Ann_R(H_{\frak a}^n(M)))H^n({\bf x}^t, M)=0,$$ for any positive integer $t$.
   Now since $H^n({\bf x}^t, M)=M/{\bf x}^tM$ we deduce that $$(\Ann_R(H_{\frak a}^0(M)))\cdots (\Ann_R(H_{\frak a}^n(M)))M\subseteq {\bf x}^tM\subseteq \frak a^tM,$$
   for any positive integer $t$, and so we have
\begin{center}
 $(\Ann_R(H_{\frak a}^0(M)))\cdots (\Ann_R(H_{\frak a}^n(M)))M\subseteq \bigcap_{t\geq 1}\frak a^tM.$
\end{center}
 Now the result follows from the  Krull's intersection theorem. \qed\\

Recall that the {\it arithmetic rank} of an ideal ${\frak a}$ in a
 Noetherian ring $R$, denoted by ${\rm ara}({\frak a})$, is the
least number of elements of $R$ required to generate an ideal which has the
same radical as ${\frak a}$, i.e.,
$${\rm ara}({\frak a})=\min\{n\in \mathbb{N}_0: \exists a_1, \dots, a_n \in R\,\,  \text{with}\,\, \Rad(a_1, \dots, a_n)=\Rad({\frak a})\}.$$
\begin{cor}
Let $R$ be a local (Noetherian) ring and $\frak a$ an ideal of $R$ with $\ara(\frak a):=n$.  Let $M$ be a finitely generated $R$-module.
Then $$(\Ann_R(H_{\frak a}^0(M)))(\Ann_R(H_{\frak a}^1(M)))\cdots (\Ann_R(H_{\frak a}^n(M)))\subseteq \Ann_RM.$$
\end{cor}
\proof The assertion follows from Proposition {\rm2.7} and the definition of $\ara(\frak a)$. \qed\\

\begin{cor}
Let $R$ be a local (Noetherian) ring, $\frak a$ an ideal of $R$ and let $M$ be a non-zero  finitely generated $R$-module
such that $\grade (\frak a, M):=g$ and $\cd (\frak a, M):=c$. Then  $$(\Ann_R(H_{\frak a}^g(M)))(\Ann_R(H_{\frak a}^{g+1}(M)))\cdots (\Ann_R(H_{\frak a}^c(M)))\subseteq \Ann_RM.$$
\end{cor}
\proof The assertion follows from Proposition {\rm2.7} and the fact that  $\cd (\frak a, M) \leq \ara(\frak a)$. \qed\\


Before bringing the next theorem  we recall the following definition.

\begin{defn}
Let $R$ be a Noetherian ring, $\frak a$ an ideal of $R$ and let $M$ be a
finitely generated $R$-module. We denote by ${\rm T}_R(\frak a, M)$ the largest submodule of $M$ such that ${\rm cd}(\frak a, {\rm T}_R(\frak a, M))<{\rm cd}(\frak a, M)$.
It follows from Lemma 2.4  that
\begin{center}
${\rm T}_R(\frak a, M)=\bigcup \{N|\,N\leq M$ and ${\rm cd} (\frak a, N)<{\rm cd} (\frak a, M)\}$.
\end{center}

\end{defn}

The following theorem improves \cite[Theorem 3.3]{Lyn}. To this end, recall that a non-zero finitely generated $R$-module $L$ is called a
relative Cohen-Macaulay module with respect to an ideal $\frak a$ of $R$ if there is precisely one non-vanishing local cohomology module of $L$ with respect to $\frak a$; that is
$\grade (\frak a, L)=\cd(\frak a, L)$.
\begin{thm}
Let $(R, \frak m)$ be a local (Noetherian)  ring and $\frak a$ an ideal of $R$. Let $M$ be a relative Cohen-Macaulay $R$-module with respect to
 $\frak a$ such that $\cd(\frak a, M):=c$. Then $$\Ann_R(H_{\frak a}^c(M))=\Ann_RM/T_R(\frak a, M)=\Ann_RM.$$
\end{thm}
\proof Since  $$\Ann_RM \subseteq \Ann_RM/T_R(\frak a, M)\subseteq \Ann_R(H_{\frak a}^c(M)),$$ it is enough for us to show that
$\Ann_R(H_{\frak a}^c(M))\subseteq \Ann_ RM.$  To do this, as $$\grade (\frak a, M)=\cd(\frak a, M),$$
 the assertion follows from Corollary {\rm2.9}. \qed\\

The next corollary improves \cite[Corollary 3.5]{Lyn}.
\begin{cor}
Let $(R, \frak m)$ be a  local (Noetherian) ring, $\frak a$ an ideal of $R$  and let $M$ be a non-zero finitely generated $R$-module.
 Suppose that $\frak a$  is generated by an $M$-regular sequence of length $n$. Then $$\Ann_R(H_{\frak a}^n(M))= \Ann_R(M/T_R(\frak a, M))=\Ann_RM.$$
\end{cor}
\proof Since $H_{\frak a}^n(M)\neq 0$ and $H_{\frak a}^i(M)= 0$ for all $i\neq n$, it follows that $\grade (\frak a, M)=\cd(\frak a, M)=n$.
Now, the assertion follows from Theorem {\rm2.11}.\qed\\
\begin{rem}
Let $R$ be a Noetherian ring and $\frak a$ an ideal of $R$. Let $M$ be a relative Cohen-Macaulay $R$-module with respect to  $\frak a$ such that $\cd(\frak a, M):=c$.
 Then, it is easy to see that  $$\Supp (H_{\frak a}^c(M))=\Supp (M)\cap V(\frak a).$$
\end{rem}
\begin{cor}
Let $(R, \frak m)$ be a  local (Noetherian) ring and $\frak a$ an ideal of $R$. Let $M$ be a non-zero finitely generated $R$-module such that ${\rm cd}(\frak a, M)=0$.
 Then $$\Ann_R(H_{\frak a}^0(M))= \Ann_R(M/T_R(\frak a, M))=\Ann_RM.$$
In particular we have $\Supp (H_{\frak a}^0(M))=\Supp (M)\subseteq V(\frak a)$.
\end{cor}
\proof The assertion follows from Theorem {\rm2.11} and Remark {\rm2.13}.\qed\\
\begin{lem}
Let $R$ be a Noetherian ring and $\frak a$ an ideal of $R$. Let $M$ be a non-zero finitely generated $R$-module such that ${\rm cd}(\frak a, M)=0$.
Then $$\Ann_R(H_{\frak a}^0(M))= \Ann_R(M/T_R(\frak a, M)).$$
\end{lem}
\proof It follows from $H_{\frak a}^0(M)\cong H_{\frak a}^0(M/T_R(\frak a, M))$  that $$ \Ann_R(M/T_R(\frak a, M))\subseteq \Ann_R(H_{\frak a}^0(M)).$$
 Now, let $x\in \Ann_R(H_{\frak a}^0(M)).$ Then from the exact sequence $$0 \longrightarrow  (0:_M xR) \longrightarrow M \stackrel{x}\longrightarrow xM \longrightarrow 0,$$
 we obtain the exact sequence $$H_{\frak a}^0(M) \stackrel{x}\longrightarrow H_{\frak a}^0(xM) \longrightarrow 0.$$
 Since $x\in \Ann_R(H_{\frak a}^0(M)),$ we get that $H_{\frak a}^0(xM)=0$, and so $xM\subseteq T_R(\frak a, M),$ as required.\qed\\
\begin{cor}
Let $(R, \frak m)$ be a  local (Noetherian) ring and $\frak a$ an ideal of $R$. Let $M$ be a non-zero finitely generated $R$-module such that ${\rm cd}(\frak a, M)=1$.
Then $$\Ann_R(H_{\frak a}^1(M))= \Ann_R(M/T_R(\frak a, M)).$$
\end{cor}
\proof Since
\begin{center}
$H_{\frak a}^1(M)\cong H_{\frak a}^1(M/H_{\frak a}^0(M))$  and $T_R(\frak a, M/H_{\frak a}^0(M))= T_R(\frak a, M)/H_{\frak a}^0(M)$,
\end{center}
 so without loss of generality we may assume that $H_{\frak a}^0(M)=0$. Now, the assertion follows  from Theorem {\rm2.11}.\qed\\
\begin{cor}
Let $(R, \frak m)$ be a local (Noetherian)  ring and $\frak a$ an ideal of $R$. Let $M$ be a relative Cohen-Macaulay $R$-module with respect to  $\frak a$ such that $\cd(\frak a, M):=c$.
 Then $$\Ass_R (R/\Ann_RM)\subseteq \{\frak p\in \Ass(M)|\,{\rm cd}(\frak a, R/\frak p)=c\}.$$
\end{cor}
\proof In view of Theorem {\rm2.11} and Lemma {\rm2.1},  there is a primary decomposition (not necessarily reduced) of $\Ann_RM,$
 \begin{eqnarray*}
\Ann_RM &=&\Ann_R(M/\cap_{\frak p_j\in \Ass_R M,\,{\rm cd}(\frak a, R/{\frak p_j})=c}N_j) \\&=& \cap_{\frak p_j\in \Ass_R M,\,{\rm cd}(\frak a, R/{\frak p_j})=c}(N_j:\,M).
\end{eqnarray*}
Therefore, there is some reduced primary decomposition of $\Ann_RM$ as $$\Ann_RM= \cap_{\frak p_j\in \Ass_R (R/\Ann_RM),\,{\rm cd}(\frak a, R/{\frak p_j})=c}(N_j:\,M).$$
 Now, the assertion follows from this.\qed\\

\section{Attached primes of local cohomology modules}

In this section we will investigate the attached prime ideals of local
cohomology modules. As the first main result, we will give a complete characterization of  the attached primes of the local cohomology module $H_{\frak a}^{\dim M-1}(M)$,
 which is a generalization of the main results of  \cite[Theorem 3.7]{ASN2} and \cite[Theorem 2.3]{He}. By using this characterization
 we show that the set $\Att_R(H_{\frak a}^{{\dim M-1}}(M))$ depends on $\Supp(M)$ only, whenever $H_{\frak a}^{\dim M}(M)=0.$
 We begin with:

\begin{defn}
 Let $L$ be an $R$-module. We say that a prime ideal  $\frak p$ of $R$ is an {\it  attached prime of}  $L$, if there exists a submodule $K$
 of $L$ such that $\frak p=\Ann_R(L/K)$ or equivalently $\frak p=\Ann_R(L/\frak p L).$ We denote by $\Att_R L$ ( resp. $\mAtt_R L$) the set of attached primes of $L$
 (resp. the set of minimal attached primes of $L$).
\end{defn}

When $M$ is {\it representable} in the sense of \cite{Ma} (e.g. Artinian or injective), our definition of $\Att_R L$ coincides with that of Macdonald \cite{Ma} and Sharp   \cite{Sh}.
The following corollary is a consequence of  Theorem 2.11.
\begin{cor}
Let $(R, \frak m)$ be a local (Noetherian)  ring and $\frak a$ an ideal of $R$. Let $M$ be a relative Cohen-Macaulay $R$-module with respect to  $\frak a$ such that $\cd(\frak a, M):=c$.
 Then $$\mAss_R M=\mAtt_R(H^c_{\frak a}(M)).$$
\end{cor}
\proof The assertion follows from Theorem 2.11 and  \cite[Lemma 3.2]{ASN2}.\qed\\

The following  theorem is our first main result of this section which extends the main results of \cite[Theorem 3.7]{ASN2} and \cite[Theorem 2.3]{He}.
  \begin{thm}
 Let $R$ be a Noetherian ring and $\frak a$ an ideal of $R$. Let $M$ be a non-zero finitely generated $R$-module of finite dimension $d$ such that $H^d_{\frak a}(M)=0$.
Then
 $$\Att_RH^{d-1}_{\frak a}(M)=\{\frak p\in \Supp (M)|\,{\rm cd}(\frak a, R/\frak p)=d-1\}.$$
\end{thm}
\proof Since $$H^{d-1}_{\frak a}(M)\cong H^{d-1}_{\frak a}(R/\Ann_RM)\otimes_{R/\Ann_RM}M,$$ so, by \cite[Lemma 2.11]{AM}, we have
\begin{eqnarray*}
\Att_{R/\Ann_RM}H^{d-1}_{\frak a}(M) &=& \Att_{R/\Ann_RM}H^{d-1}_{\frak a}(R/\Ann_RM)\cap \Supp_{R/\Ann_RM}(M) \\&=& \Att_{R/\Ann_RM}H^{d-1}_{\frak a}(R/\Ann_RM).
\end{eqnarray*}
So, in view of \cite[Theorem 3.7]{ASN2}, we have $$\Att_{R/\Ann_RM}H^{d-1}_{\frak a}(M)=\{\frak p/\Ann_RM\in \Spec (R/\Ann_RM)|\,{\rm cd}(\frak a, R/\frak p)= d-1\}.$$
 Now, as $$\Att_RH^{d-1}_{\frak a}(M)=\{\frak p\in \Supp M|\,\frak p/\Ann_RM\in\Att_{R/\Ann_RM}H^{d-1}_{\frak a}(M)\},$$ the assertion follows.\qed\\
\begin{cor}
Let $R$ be a Noetherian ring and $\frak a$ an ideal of $R$. Let $M$ and $N$ be two non-zero finitely generated $R$-modules and $\dim M=d$ and $H_{\frak a}^d(M)=0.$
If $\,\Supp (N) \subseteq \Supp (M)$, then $\Att_R(H_{\frak a}^{d-1}(N)) \subseteq \Att_R(H_{\frak a}^{d-1}(M))$.
\end{cor}
\proof By Lemma {\rm2.4}, we can (and do) assume that ${\rm cd}(\frak a, N)= {\rm cd}(\frak a, M)=d-1.$ Now, the assertion follows from Theorem {\rm3.3}.\qed\\

The next consequence of Theorem 3.3 shows that the attached primes of the local cohomology module $H_{\frak a}^{\dim M-1}(M)$,
depends only on $\Supp(M)$, whenever $H_{\frak a}^{\dim M}(M)=0.$
\begin{cor}
Let $R$ be a Noetherian ring and $\frak a$ an ideal of $R$. Let $M$ and $N$ be two  finitely generated $R$-modules with $\dim M=d$ and
 $H_{\frak a}^d(M)=0.$ If $\Supp (M)=\Supp (N)$, then $\Att_R(H_{\frak a}^{d-1}(N))=\Att_R(H_{\frak a}^{d-1}(M)).$
\end{cor}
\proof The assertion follows from Lemma {\rm2.4} and Corollary {\rm3.4}.\qed\\

\begin{lem}
Let $R$ be a Noetherian domain and $\frak a$ an ideal of $R$ such that ${\rm cd}(\frak a, R)=1$. Then $\Ann_RH^1_{\frak a}(R)=0$.
\end{lem}
\proof
Suppose in contrary that $x(\neq0)\in \Ann_RH^1_{\frak a}(R).$ Then, Since $H^0_{\frak a}(R)=0,$ by the exact sequence
 $$0 \longrightarrow R \stackrel{x}\longrightarrow R \longrightarrow R/xR \longrightarrow 0,$$ we get $H^0_{\frak a}(R/xR)\cong H^1_{\frak a}(R).$
 This is a contradiction, because, $H^0_{\frak a}(R/xR)$ is a finitely generated $R$-module, but $H^1_{\frak a}(R)$, is not.\qed\\
\begin{prop}
Let $R$ be a Noetherian ring and $\frak a$ an ideal of $R$ such that ${\rm cd}(\frak a, R)=1$. Then

$$\Att_RH^1_{\frak a}(R)= \{\frak p\in \Spec R|\,{\rm cd}(\frak a, R/\frak p)=1\}.$$
\end{prop}
\proof In view of \cite[Theorem 3.3]{ASN2}, we have $$\Att_RH^1_{\frak a}(R)\subseteq \{\frak p\in \Spec R|\,{\rm cd}(\frak a, R/\frak p)=1\}.$$
 Now, let $\frak p$ be a prime ideal of $R$ such that ${\rm cd}(\frak a, R/\frak p)=1.$ Then, by Lemma {\rm3.7}, $\Ann_{R/\frak p}H^1_{\frak a}(R/\frak p)=0.$
 So, we have $\Ann_RH^1_{\frak a}(R/\frak p)=\frak p.$ Thus by Definition {\rm3.1}, we have $\frak p\in\Att_RH^1_{\frak a}(R/\frak p).$ Now, from the exact sequence
 $$0 \longrightarrow \frak p \longrightarrow R \longrightarrow R/{\frak p} \longrightarrow 0,$$
 and the right exactness of  $H^1_{\frak a}(\cdot)$, we deduce that $\frak p \in \Att_RH^1_{\frak a}(R)$, as required.\qed\\

 The following result gives a partial answer to \cite[Question (i)]{ASN2}, in the case  $\cd(\frak a, M)=1$.
 \begin{thm}
 Let $R$ be a Noetherian ring and $\frak a$ an ideal of $R$. Let $M$ be a non-zero finitely generated $R$-module such that ${\rm cd}(\frak a, M)=1$.
Then
 $$\Att_RH^1_{\frak a}(M)=\{\frak p\in \Supp (M)|\,{\rm cd}(\frak a, R/\frak p)=1\}.$$
\end{thm}
\proof The proof is similar to the proof of Theorem {\rm3.3}, by using Proposition {\rm3.7} instead of \cite[Theorem 3.7]{ASN2}.\qed\\

The next  corollary reproves \cite[Corollary 2.4]{He}.
\begin{cor}
Let $R$ be a local (Noetherian) ring and $\frak a$ an ideal of $R$. Let $M$ be a non-zero finitely generated $R$-module such that ${\rm cd}(\frak a, M)\leq1$.
Then $$\Att_RH^1_{\frak a}(M)=\Supp (M) \setminus V(\frak a).$$ In particular, if $x\in R$, then $\Att_RH^1_{xR}(M)=\Supp (M) \setminus V(xR).$
\end{cor}
\proof The assertion follows from Corollary {\rm2.14} and Theorem {\rm3.8}.\qed\\
\begin{rem}
Let $R$ be a Noetherian ring and $\frak a$ an ideal of $R$. Let $M$ be a non-zero finitely generated $R$-module such that ${\rm cd}(\frak a, M)=0$.
Then it is not hard to see that  $\frak p\in V(\frak a)$ if and only if ${\rm cd}(\frak a, R/\frak p)=0$, for each $\frak p\in \Supp(M)$.
\end{rem}
\begin{lem}
Let $R$ be a Noetherian ring and $\frak a$ an ideal of $R$. Let $M$
be a non-zero finitely generated $R$-module such that ${\rm
cd}(\frak a, M)=0$. Then  $$\Att_RH^0_{\frak a}(M)=\{\frak p\in
\Supp (M)|\,{\rm cd}(\frak a, R/\frak p)=0\}.$$
\end{lem}
\proof Since $H^0_{\frak a}(M)$ is a finitely generated $R$-module, we have $$\Att_RH^0_{\frak a}(M)=\Supp (H^0_{\frak a}(M)).$$
Now,  the assertion follows from of Remarks {\rm2.13} and {\rm3.10}.\qed\\

 The following theorem shows that  \cite[Question (i)]{ASN2}, is true in the case  $\dim M\leq 3$.
\begin{thm}
 Let $R$ be a Noetherian ring and $\frak a$ an ideal of $R$. Let $M$ be a non-zero finitely generated $R$-module with $\dim M\leq 3$ and $c:={\rm cd}(\frak a, M).$
Then
 $$\Att_RH^c_{\frak a}(M)=\{\frak p\in \Supp (M)|\,{\rm cd}(\frak a, R/\frak p)=c\}.$$
\end{thm}
\proof The result follows from Lemma {\rm3.11}, Theorems {\rm3.3}, {\rm3.8} and \cite[Corollary 3.4]{ASN2}.\qed\\

It is natural to ask about the attached primes of $\Att_RH^{\cd(\frak a, M)}_{\frak a}(M)$ for arbitrary $R$-module $M$ (not necessarily finitely generated).
The following is our second main result of this section which gives a partial answer to this question. As a consequence we extend the main results of   \cite{DY} and \cite{Ru}.

\begin{thm}
Let $R$ be a Noetherian ring and $\frak a$ an ideal of $R$. Let $M$ be a non-zero $R$-module (not necessarily finitely generated) such that
$c:=\cd(\frak a, M)={\rm cd}(\frak a, R/\Ann_RM)$  is finite.
Then
 $$\Att_RH^c_{\frak a}(M)\subseteq\{\frak p\in \V(\Ann_RM)|\,{\rm cd}(\frak a, R/\frak p)=c\}.$$
\end{thm}
\proof Let $\frak p\in\Att_RH^c_{\frak a}(M)$. Then clearly we have $\Ann_RM\subseteq\frak p$ and so
$$\Supp(R/\frak p)\subseteq \Supp(R/\Ann_RM).$$
Hence  by using  Lemma {\rm2.4} we see that ${\rm cd}(\frak a, R/\frak p)\leq c$. On the other hand,  it follows from $\frak p\in \Att_RH^c_{\frak a}(M)$
 that $$\frak p/\Ann_RM\in \Att_{R/\Ann_RM}H^c_{\frak a}(M),$$ and so
$$H^c_{\frak a}(M)\otimes_{R/\Ann_RM}R/\frak p\neq 0.$$ Now, as $$H^c_{\frak a}(M)\cong H^c_{\frak a}(R/\Ann_RM)\otimes_{R/\Ann_RM}M,$$ it follows that
$$H^c_{\frak a}(R/\Ann_RM)\otimes_{R/\Ann_RM}M\otimes_{R/\Ann_RM}R/\frak p\neq 0.$$
Consequently $$H^c_{\frak a}(R/\frak p)\otimes_{R/\Ann_RM}M\neq 0,$$  and thus $ H^c_{\frak a}(R/\frak p)\neq 0$, as required.\qed\\

The final result extends the main result of  \cite[Theorem B]{DY}.
\begin{cor}
Let $R$ be an arbitrary Noetherian ring and $\frak a$ an ideal of $R$. Let $M$ be a non-zero $R$-module (not necessarily finitely generated) such that $d:=\dim M=\dim R$  is finite.
Then
 $$\Att_RH^d_{\frak a}(M)\subseteq\{\frak p\in \Ass_R M|\,{\rm cd}(\frak a, R/\frak p)=d\}.$$
\end{cor}
\proof We may assume that $\cd(\frak a , M)=d$. Then by Theorem {\rm 3.13} it is enough to show that $\cd(\frak a , R/\Ann_RM)=d$. To do this, as $\Supp (M)\subseteq\Supp(R/\Ann_RM)$,
it follows from \cite[Theorem 1.4]{DY1} that  $$d=\cd(\frak a , M)\leq\cd(\frak a , R/\Ann_RM)\leq\cd(\frak a , R)\leq\dim R=d,$$ as required.
   \qed\\

\begin{center}
{\bf Acknowledgments}
\end{center}
The authors are deeply grateful to the referee for  careful reading and many useful suggestions. We would like to thank  Prof. Kamran  Divaani-Aazar for useful discussions,
 and also the Azarbaijan Shahid Madani University, for the financial support.

\end{document}